\documentclass[12pt]{article}

\usepackage{amsfonts}

\newtheorem{thm}{Theorem}[section]
\newtheorem{lemma}[thm]{Lemma}
\newtheorem{cor}[thm]{Corollary}

\newtheorem{prop}[thm]{Proposition}
\newtheorem{conjecture}{Conjecture}

\newcommand{\be}{\begin{equation}}
\newcommand{\ee}{\end{equation}}
\newcommand{\openbox}{\leavevmode
  \hbox to8pt{\hfil\vrule\vbox to6pt{\hrule width6pt\vfil\hrule}\vrule}}

\newcommand{\qed}{\hbox to5pt{ } \hfill \openbox\bigskip\medskip}

\newcommand{\Zp}{\mathbb Z _p}

\newcommand{\cF}{\mbox{$\cal F$}}

\newcommand{\ve}[1]{\mathbf{#1}}

\newcommand{\R}{\mathbb R}

\newcommand{\F}{\mathbb F}

\title{Non-uniform skew versions of Bollob\'as' Theorem}
\author{G\'abor Heged\"{u}s
\\{\normalsize  \'Obuda University}
\\{\normalsize B\'ecsi \'ut 96/B, Budapest, Hungary, H-1032}
\\{\normalsize hegedus.gabor@uni-obuda.hu}
}

\begin{document}
\maketitle

\begin{abstract}
Let $A_1, \ldots ,A_m$ and $B_1, \ldots ,B_m$ be subsets of $[n]$ and let $t$ be a non-negative integer with the following property:  $|A_i \cap B_i|\leq t$ for each $i$ and 
$|A_i\cap B_j|>t$ whenever $i< j$. Then
$m\leq  2^{n-t}$.  Our proof uses Lov\'asz' tensor product method. 

We prove the following  skew version of Bollob\'as' Theorem. 
Let $A_1, \ldots ,A_m$ and  $B_1, \ldots ,B_m$ be finite sets of $[n]$ satisfying the conditions $A_i \cap B_i =\emptyset$ for each $i$ and 
$A_i\cap B_j\ne \emptyset$ for each $i< j$. 
Then 
$$
\sum_{i=1}^m \frac{1}{{|A_i|+|B_i| \choose |A_i|}}\leq n+1.
$$
Both upper bounds are sharp. 
\end{abstract}

\section{Introduction}

Bollob\'as proved in \cite{B} the following two remarkable results in extremal combinatorics.
\begin{thm} \label{Boll}
Let $A_1, \ldots ,A_m$ and  $B_1, \ldots ,B_m$ be finite sets satisfying the conditions 
\begin{itemize}
\item[(i)] $A_i \cap B_i =\emptyset$ for each $1\leq i \leq m$;
\item[(ii)] $A_i\cap B_j\ne \emptyset$ for each $i\ne j$ ($1\leq i, j \leq m$).
\end{itemize}
Then 
$$
\sum_{i=1}^m \frac{1}{{|A_i|+|B_i| \choose |A_i|}}\leq 1.
$$
\end{thm}

\begin{thm} \label{Boll2}
Let $A_1, \ldots ,A_m$ be $r$-element sets and  $B_1, \ldots ,B_m$ be $s$-element sets such that 
\begin{itemize}
\item[(i)] $A_i \cap B_i =\emptyset$ for each $1\leq i \leq m$;
\item[(ii)] $A_i\cap B_j\ne \emptyset$ for each $i\ne j$ ($1\leq i, j \leq m$).
\end{itemize}
Then
$$
m\leq {r+s \choose s}.
$$
\end{thm}

Lov\'asz used in \cite{L1} tensor product methods to prove the following skew version of Bollob\'as' Theorem for subspaces.

\begin{thm} \label{Lovasz}
Let $\F$ be an arbitrary field.  Let  $U_1, \ldots ,U_m$ be $r$-dimensional and $V_1, \ldots ,V_m$ be $s$-dimensional subspaces of a linear  space $W$ over the field $\F$. Assume that 
\begin{itemize}
\item[(i)] $U_i \cap V_i =\{0\}$ for each $1\leq i \leq m$;
\item[(ii)] $U_i\cap V_j\ne \{0\}$ whenever $i< j$ ($1\leq i, j \leq m$).
\end{itemize}
Then
$$
m\leq {r+s \choose r}.
$$
\end{thm}

F\"uredi generalized Lov\'asz' result in \cite{F} and obtained the following threshold variant of Theorem  \ref{Lovasz}. 

\begin{thm} \label{Furedi1}
Let $\F$ be an arbitrary field.  Let  $U_1, \ldots ,U_m$ be $r$-dimensional and $V_1, \ldots ,V_m$ be $s$-dimensional subspaces of a linear  space $W$ over the field $\F$. Assume that for some $t\geq 0$
\begin{itemize}
\item[(i)] $\dim(U_i \cap V_i)\leq t$ for each $1\leq i \leq m$;
\item[(ii)] $\dim(U_i\cap V_j)>t$ whenever $i< j$ ($1\leq i, j \leq m$).
\end{itemize}
Then
$$
m\leq {r+s-2t \choose r-t}.
$$
\end{thm}

F\"uredi derived the following combinatorial result from this subspace version of Bollob\'as' Theorem.

\begin{thm} \label{Furedi2}
Let $A_1, \ldots ,A_m$ be $r$-element sets and  $B_1, \ldots ,B_m$ be $s$-element sets. Assume that for some $t\geq 0$
\begin{itemize}
\item[(i)] $|A_i \cap B_i|\leq t$ for each $1\leq i \leq m$;
\item[(ii)] $|A_i\cap B_j|>t$ whenever $i< j$ ($1\leq i, j \leq m$).
\end{itemize}
Then
$$
m\leq {r+s-2t \choose r-t}.
$$
\end{thm}

In this paper one of our main aim is to give non-uniform versions of F\"uredi's results.

First we state a non-uniform variant of Lov\'asz' result about subspaces.
\begin{thm} \label{main}
Let $\F$ be an arbitrary field.  Let  $U_1, \ldots ,U_m$ and $V_1, \ldots ,V_m$ be  subspaces of an $n$-dimensional linear  space $W$ over the field $\F$. Assume that 
\begin{itemize}
\item[(i)] $U_i \cap V_i =\{0\}$ for each $1\leq i \leq m$;
\item[(ii)] $U_i\cap V_j\ne \{0\}$ whenever $i< j$ ($1\leq i, j \leq m$).
\end{itemize}
Then
$$
m\leq 2^n.
$$
\end{thm}

Next we give  a  non-uniform variant of  F\"uredi's  Theorem  \ref{Furedi1}. Our proof based on the tensor product method.

\begin{thm} \label{main2}
Let $\F$ be an arbitrary field.  Let  $U_1, \ldots ,U_m$ and $V_1, \ldots ,V_m$ be  subspaces of an $n$-dimensional linear  space $W$ over the field $\F$. Assume that for some $0\leq t\leq n$
\begin{itemize}
\item[(i)] $\dim(U_i \cap V_i)\leq t$ for each $1\leq i \leq m$;
\item[(ii)] $\dim(U_i\cap V_j)>t$ whenever $i< j$ ($1\leq i, j \leq m$).
\end{itemize}
Then
$$
m\leq  2^{n-t}.
$$
\end{thm}

The following Corollary is immediate.

\begin{cor} \label{main22}
Let $A_1, \ldots ,A_m$ and $B_1, \ldots ,B_m$ be subsets of $[n]$. Assume that for some  $0\leq t\leq n$
\begin{itemize}
\item[(i)] $|A_i \cap B_i|\leq t$ for each $1\leq i \leq m$;
\item[(ii)] $|A_i\cap B_j|>t$ whenever $i< j$ ($1\leq i, j \leq m$).
\end{itemize}
Then
$$
m\leq  2^{n-t}.
$$
\end{cor}

The following example shows that Corollary \ref{main22} is sharp. 
Let $T:=\{n-t+1,\ldots ,n\}$.  First list all pairs $(A,[n-t]\setminus A)$ with $A\subseteq [n-t]$, where we sort them with decreasing cardinality of the first element $A$. Then consider the list $(A\cup T, [n]\setminus A)$ with the same order. 

Finally we prove the following  skew version of Bollob\'as' Theorem.

\begin{thm} \label{main3}
Let $A_1, \ldots ,A_m$ and  $B_1, \ldots ,B_m$ be finite sets of $[n]$ satisfying the conditions 
\begin{itemize}
\item[(i)] $A_i \cap B_i =\emptyset$ for each $1\leq i \leq m$;
\item[(ii)] $A_i\cap B_j\ne \emptyset$ for each $i< j$ ($1\leq i, j \leq m$).
\end{itemize}
Then 
$$
\sum_{i=1}^m \frac{1}{{|A_i|+|B_i| \choose |A_i|}}\leq n+1.
$$
\end{thm}

Babai and Frankl gave the following example in \cite{BF} Exercise 5.1.1. , which shows that Theorem \ref{main3} is sharp. This example is a list of all pairs $(A,[n]\setminus A)$ with $A\subseteq [n]$, sorted them with decreasing cardinality of the first element. Let $a_i:=|A_i|$ for each $1\leq i\leq 2^n$. Then
$$                    
\sum_{i=1}^{2^n} \frac{1}{{n\choose a_i}}=\sum_{j=0}^n \frac{{n\choose j}}{{n\choose j}}=n+1.
$$
 
The proof of Theorem \ref{main3} is completely elementary. We modified slightly one of the proof of Bollob\'as' Theorem. 

\section{Preliminaries}

The following Lemma was proved in \cite{FT} Lemma 26.14. 
\begin{lemma} \label{subs_gen_pos}
Let $n\geq k$ and $t=n-k$. Let $\F$ be a field. Let $V$ be an $n$-dimensional vector space. Let $W_1, \ldots ,W_m$ be subspaces of $V$ with $\dim(W_i)<n$ for each $i$. Then there exists a $k$-dimensional subspace $V'$ such that 
$$
\dim(W_i\cap V')=\mbox{\rm max}(\dim(W_i)-t,0)
$$ 
for each $1\leq i\leq m$.          
\end{lemma}

This subspace $V'$ guaranteed by Lemma  \ref{subs_gen_pos} is called to be in {\em general position} with respect to the subspaces $W_1, \ldots ,W_m$.

The following Proposition appears as a Triangular Criterion in \cite{BF} Proposition 2.9.
  
\begin{prop} \label{triang}
Let $\F$ be an arbitrary field. Let $W,T$ be linear spaces over $\F$ and $\Omega$ an arbitrary set. Let $f:W\times \Omega\to T$ be a function  which is linear in the first variable. For $i=1,\ldots ,m$, let $\ve w_i\in W$ and $a_i\in \Omega$ be such that 
\[ f(\ve w_i,a_j)\left\{ \begin{array}{ll}
\neq 0, & \textrm{if $i=j$;} \\
=0, & \textrm{if $i<j$.}
\end{array} \right. \]
Then the vectors $\ve w_1, \ldots ,\ve w_n$ are linearly independent.
\end{prop}

Let $\F$ be a field. Let $V$ be an $n$-dimensional vector space over $\F$. Write $\bigwedge V$ for the exterior algebra of $V$. 

Let $E=\{\ve e_1, \ldots ,\ve e_n\}$ denote the standard basis for $V$.

For $A\subseteq [n]$, write $f_A:= \wedge_{i\in A} \ve e_i \in \bigwedge V$. Let $F:= \{f_A:~ A\subseteq [n]\}$. It is a  well-known basic fact (see \cite{HH} Chapter 5) that  $F\subseteq \bigwedge V$ is a basis of the exterior algebra, hence $\dim(\bigwedge V)=2^n$. 

For a $k$-dimensional subspace $T\leq V$ define $\wedge T\in \bigwedge V$ by selecting a basis $\ve t_1, \ldots ,\ve t_k$ of $T$ and setting
$$
\wedge T:=\ve t_1 \wedge \ldots \wedge \ve t_k.
$$
Clearly this definition depends on the special choice of the basis   $\ve t_1, \ldots ,\ve t_k$, but it is a well-known statement that 
it can only vary by a nonzero scalar factor.

The next result links the wedge products to intersection conditions.

\begin{lemma} \label{wedge}
Let $\F$ be a field. Let $W$ be an $n$-dimensional vector space over $\F$.
Let $U$ and $V$ be subspaces of $W$. Then $(\wedge U)\wedge (\wedge V)=0$ iff $U\cap V\neq 0$.
\end{lemma}
\section{Proofs}

{\bf Proof of Theorem \ref{main}:}

Let $u_i:=\wedge U_i$ and  $v_i:=\wedge V_i$ for each $i$.
If we combine Lemma \ref{wedge} with the conditions of Theorem  \ref{main}, we get that
\[ u_i\wedge v_j\left\{ \begin{array}{ll}
\neq 0, & \textrm{if $i=j$;} \\
=0, & \textrm{if $i<j$.}
\end{array} \right. \]
It follows from Proposition \ref{triang} that the vectors $u_1, \ldots , u_m$ are linearly independent. This means that $m\leq 2^n$. \qed

{\bf Proof of Theorem \ref{main2}:}
Let $W_0$ be a subspace of codimension $t$ in general position with respect to the subspaces $U_i$, $V_i$ and $U_i\cap V_i$ (here $1\leq i\leq m$). It follows from Lemma \ref{subs_gen_pos} that there exists such $W_0$ subspace  in general position. 

Then
\begin{itemize}
\item[(i)] $\dim(U_i \cap V_i\cap W_0)= 0$ for each $1\leq i \leq m$;
\item[(ii)] $\dim(U_i\cap V_j\cap W_0)>0$ whenever $i< j$ ($1\leq i, j \leq m$).
\end{itemize}
Observe that the conditions (i) and (ii) guarantee that for the subspaces $U_i\cap W_0$, $V_i\cap W_0$ and $W_0$ the conditions of Theorem \ref{main} hold, but now  with $n-t$ in the role of $n$. Consequently $m\leq  2^{n-t}$. \qed

{\bf Proof of Corollary \ref{main22}:}

Let $X$ be a finite subset containing all the $A_i$ and $B_j$. Suppose that $|X|=n$ and let $W:={\R}^n$. Associate with each $x\in X$ a member $\ve e_x$ of a fixed basis of the linear space $W$.  Associate with each subset $T\subseteq X$ the subspace $W(T):=\mbox{\rm span}\{\ve e_x:~ x\in T\}$. 

Then the subspaces $W(A_i)$  and $W(B_i)$ satisfy the conditions of Theorem \ref{main2}, consequently we get the corresponding bound on $m$. \qed

{\bf Proof of Theorem \ref{main3}:}

Let $X$ be the union of all sets $A_i\cup B_j$. We prove by induction on $n$, where $n=|X|$. 

If $n=1$, then it is easy to verify the statement. 
Assume it holds for $n-1$ and we prove it for $n$. Let $a_i:=|A_i|$ and $b_i:=|B_i|$ for each $i$.

It is enough to prove that
$$
\sum_{i:~ A_i\neq X} \frac{1}{{a_i+b_i \choose a_i}} \leq n,
$$
because clearly 
$$
\sum_{i:~ A_i= X} \frac{1}{{a_i+b_i \choose a_i}} =1.
$$

Let $x\in X$ be a fixed element and define the family of pairs
$$
{\cF}_x:=\{(A_i,B_i\setminus \{x\}):~x\notin A_i\}.
$$
Here we can consider this family ${\cF}_x$ as an ordered family, where the order on this family is inherited from the order on the original family $(A_i,B_i)_{1\leq i\leq m}$.

Clearly the base set of each of these families ${\cF}_x$ has less than $n$ elements, consequently we can apply our induction hypothesis for each of them and we can sum the corresponding inequalities.
We get that
$$
\sum_{x\in X} \sum_{(A',B')\in {\cF}_x}  \frac{1}{{a'+b' \choose a'}}\leq \sum_{x\in X} n=n^2,
$$
where $a'=|A'|$ and $b'=|B'|$.

It is easy to see that the previous sum counts $n-a_i-b_i$ times the term $1/{a_i+b_i \choose a_i}$ corresponding to points $x\notin A_i\cup B_i$ and $b_i$ times the term $1/{a_i+b_i-1 \choose a_i}$ corresponding to points $x\in B_i$. Consequently we get that
$$
\sum_{x\in X} \sum_{(A',B')\in {\cF}_x}  \frac{1}{{a'+b' \choose a'}}=\sum_{i:~ A_i\neq X}\frac{ n-a_i-b_i}{{a_i+b_i \choose a_i}} + \frac{ b_i}{{a_i+b_i-1 \choose a_i}}=
$$
$$
=\sum_{i:~ A_i\neq X} \frac{ n-a_i-b_i}{{a_i+b_i \choose a_i}} + \frac{ b_i}{\frac{b_i}{a_i+b_i}\cdot{a_i+b_i \choose a_i}}=
$$
$$
=n\sum_{i:~ A_i\neq X} \frac{1}{{a_i+b_i \choose a_i}}.
$$
Hence we obtain that
$$
n\sum_{i:~ A_i\neq X} \frac{1}{{a_i+b_i \choose a_i}}\leq n^2.
$$
Dividing both side by $n$ we get our result. \qed

\section{Concluding remarks}

We can raise  the following natural question:  Theorem \ref{main3} is sharp and we described an optimal example. Determine the extremal examples for Theorem \ref{main3}, i.e., all such $A_1, \ldots ,A_m$ and $B_1, \ldots ,B_m$ subsets of $[n]$ satisfying the conditions 
\begin{itemize}
\item[(i)] $A_i \cap B_i =\emptyset$ for each $1\leq i \leq m$;
\item[(ii)] $A_i\cap B_j\ne \emptyset$ for each $i< j$ ($1\leq i, j \leq m$);
\item[(iii)] $\sum_{i=1}^m \frac{1}{{|A_i|+|B_i| \choose |A_i|}}= n+1$.
\end{itemize}

Here we conjecture the following.

\begin{conjecture} \label{Hconj}
Let $A_1, \ldots ,A_m$ and $B_1, \ldots ,B_m$ be subsets of $[n]$ satisfying the conditions 
\begin{itemize}
\item[(i)] $A_i \cap B_i =\emptyset$ for each $1\leq i \leq m$;
\item[(ii)] $A_i\cap B_j\ne \emptyset$ for each $i< j$ ($1\leq i, j \leq m$);
\item[(iii)] $\sum_{i=1}^m \frac{1}{{|A_i|+|B_i| \choose |A_i|}}= n+1$.
\end{itemize}
Then this family of pairs $(A_i,B_i)_{1\leq i\leq m}$ is a list of all pairs $(A,[n]\setminus A)$ with $A\subseteq [n]$, sorted them with decreasing cardinality of the first element.
\end{conjecture}

\end{document}